\documentclass[11pt]{article}
\usepackage{amsmath, amsfonts, amssymb, amsthm, newlfont, graphicx}

\usepackage{amssymb,latexsym}
\usepackage[mathscr]{eucal}
\usepackage{setspace}
\usepackage{graphics}
\usepackage{array}

\begin{document}
\title{Integrable Spatiotemporally Varying NLS, PT-Symmetric NLS, and DNLS Equations: Generalized Lax Pairs and Lie Algebras}
\author{Matthew Russo and S. Roy Choudhury\\  \small  Department of Mathematics, University of Central Florida, Orlando, FL  32816-1364 USA\\
Corresponding author email: choudhur@cs.ucf.edu }        
\date{\today}       
\maketitle

\centerline{ \noindent\Large\textbf{Abstract}}
This paper develops two approaches to Lax-integrbale systems with spatiotemporally varying coefficients.
A technique based on extended Lax Pairs is first considered to derive variable-coefficient generalizations of various Lax-integrable NLPDE hierarchies recently introduced in the literature. As illustrative examples, we consider generalizations of the NLS and DNLS equations, as well as
a PT-symmetric version of the NLS equation. It is demonstrated that the techniques yield Lax- or S-integrable NLPDEs with both time- AND space-dependent
coefficients which are thus more general than almost all cases considered earlier via other methods such as the
Painlev\'e Test, Bell Polynomials, and various similarity methods.

However, this technique, although operationally effective, has the significant disadvantage that, for any integrable system with spatiotemporally varying coefficients, one must 'guess' a generalization of the structure of the known Lax Pair for the corresponding system with constant coefficients. Motivated by the somewhat arbitrary nature of the above procedure, we therefore next attempt to systematize the derivation of Lax-integrable sytems with variable coefficients. An ideal approach would be a method which does not require knowledge of the Lax pair to an associated constant coefficient system, and also involves little to no guesswork. Hence we attempt to apply the Estabrook-Wahlquist (EW) prolongation technique, a relatively self-consistent procedure requiring little prior infomation. However, this immediately requires that the technique be significantly generalized or broadened in several different ways, including
solving matrix partial differential equations instead of algebraic ones as the structure of the Lax Pair is deduced systematically following the standard Lie-algebraic procedure of proceeding downwards from the coefficient of the highest derivative. The same is true while finding the explicit forms for the various 'coefficient' matrices which occur in the procedure, and which must satisfy the various constraint equations which result at various stages of the calculation. 

The new and extended EW technique which results is illustrated by algorithmically deriving generalized Lax-integrable versions of NLS,
PT-symmetric NLS, and DNLS equations.
\vspace{0.25in}

Key Words: Generalizing Lax or S-integrable equations, spatially and temporally-dependent coefficients, generalized Lax Pairs, extended Estabrook-Wahlquist method.

\section{Introduction}

     Variable Coefficient Korteweg de Vries (vcKdV) and Modified Korteweg de Vries (vcMKdV)equations have a long history dating from their derivation in various applications\cite{K1}-\cite{K10}. However, almost all studies, including those which derived exact solutions by a variety of techniques, as well as those which considered integrable sub-cases and various integrability properties by methods such as Painlev\'e analysis, Hirota's method, and Bell Polynomials treat vcKdV equations with coefficients which are functions of the time only. For instance, for generalized variable coefficient NLS (vcNLS) equations, a particular coefficient is usually taken to be a function of $x$ \cite{K11}, as has also been sometimes done for vcMKdV equations\cite{K12}. The papers \cite{K13}-\cite{Khawaja} are somewhat of an exception in that they treat vcNLS equations having coefficients with general $x$ and $t$ dependences. Variational principles, solutions, and other integrability properties have also been considered for some of the above variable coefficient NLPDEs in cases with time-dependent coefficients.

In applications, the coefficients of vcKdV or vcNLS equations may include spatial dependence, in addition to the temporal variations that have been extensively considered using a variety of techniques. Both for this reason, as well as
for their general mathematical interest, extending integrable hierarchies of nonlinear PDEs (NLPDEs) to include {\it both} spatial and temporal dependence of the coefficients is worthwhile.

Given the above, we compare two methods for deriving the integrability conditions of both a general form of variable-coefficient MKdV (vcMKdV) equation, as well as a general, variable-coefficient KdV (vcKdV) equation. In all of these cases,the coefficients are allowed to vary in space AND time.

The first method employed here is based on directly establishing Lax integrability (or S-integrability to use the technical term) as detailed in the following sections. As such, it is fairly 
general, although subject to the ensuing equations being solvable. We should stress that the computer algebra involved is quite challenging, and an order of magnitude beyond that encountered for integrable, constant coefficient NLPDEs.

However, this first technique, although operationally effective, has the significant disadvantage that, for any integrable system with spatiotemporally varying coefficients, one must 'guess' a generalization of the structure of the known Lax Pair for the corresponding system with constant coefficients. This involves
replacing constants in the Lax Pair for the constant coefficient integrable system, including powers of the spectral parameter, by functions. Provided that one has guessed correctly and generalized the constant coefficient system's Lax Pair sufficiently, and this is  of course hard to be sure of 'a priori', one may then
proceed to systematically deduce the Lax Pair for the corresponding variable-coefficient integrable system.

Motivated by the somewhat arbitrary nature of the above procedure, we next attempt to systematize the derivation of Lax-integrable sytems with variable coefficients. Of the many techniques which have been employed for constant coefficient integrable systems, the Estabrook-Wahlquist (EW) prolongation technique \cite{EW1}-\cite{EW4}is among the most self-contained. The method directly proceeds to attempt construction of the Lax Pair or linear spectral problem, whose compatibility condition is the integrable system under discussion. While not at all guaranteed to work, any successful implementation of the technique means that Lax-integrability has already been verified during the procedure, and in addition the Lax Pair is algorithmically obtained. If the technique fails, that does not necessarily imply non-integrability  of the equation contained in the compatibility condition of the assumed Lax Pair. It may merely mean that some of the starting assumptions may not be appropriate or general enough.

Hence we attempt to apply the Estabrook-Wahlquist (EW) method as a second, more algorithmic, technique to generate a variety of such integrable systems with such spatiotemporally varying coefficients. However, this immediately requires that the technique be significantly generalized or broadened in several different ways which we then develop and outline, before illustrating this new and extended method with examples.

The outline of this paper is as follows. In Section 2, we briefly review the Lax Pair method and its modifications for variable-coefficient NLPDEs, and then apply it to the vcNLS and vc Pt-symmetric NLS systems. In section 3 we consider an analogous treatment of a generalized vcDNLS equation.  In section 4, we lay out the extensions required to apply the EW procedure to Lax-integrable systems
with spatiotemporally varying coefficients. Sections 5 and 6 then illustrate this new, extended EW method in detail for Lax-integrable versions of the 
NLS/PTNLS and generalized DNLS systems respectively, each with spatiotemporally varying coefficients. We also illustrate that this generalized
EW procedure algorithmically generates the same results as those obtained in a more ad hoc manner in Sections 2 and 3.
 Section 7 briefly reviews the results and conclusions, and directions for possible future work.

More involved algebraic details, which are integral to both the procedures employed here, are relegated to the appendices.

\section{Extended Lax Pair method and application to the vc-NLS and vcPT-Symmetric NLS Systems}

In the Lax pair method \cite{K15} - \cite{K16} for solving and determining the integrability conditions for nonlinear partial differential equations (NLPDEs) a pair of $n\times n$ matrices, $\textbf{U}$ and $\textbf{V}$ needs to be derived or constructed. The key component of this construction is that the integrable nonlinear PDE under consideration must be contained in, or result from, the compatibility of the following two linear Lax equations (the Lax Pair)
\begin{eqnarray}
\Phi_{x} &=& U\Phi \\
\Phi_{t} &=& V\Phi
\end{eqnarray}
where $\Phi$ is an eigenfunction, and $\textbf{U}$ and $\textbf{V}$ are the time-evolution and spatial-evolution (eigenvalue problem) matrices. 

From the cross-derivative condition (i.e. $\Phi_{xt} = \Phi_{tx}$) we get
\begin{equation} \label{ZCC1}
U_{t}-V_{x}+[U,V] = \dot{0}
\end{equation}
known as the zero-curvature condition where $\dot{0}$ is contingent on $q(x,t)$ and $r(x,t)$ being solutions to the system of nonlinear PDEs. A Darboux transformation can then be applied to the linear system to obtain solutions from known solutions and other integrability properties of the integrable NLPDE. 

We first consider the variable coefficient cubic nonlinear Schrodinger and PT-symmetric nonlinear Schrodinger equations given by

\begin{equation}
iq_{t}(x,t) + f(x,t)q_{xx}(x,t) + g(x,t)|q(x,t)|^{2}q(x,t) + v(x,t)q(x,t) + i\gamma(x,t)q(x,t) = 0
\end{equation}

\noindent
and

\begin{equation}
iq_{t}(x,t) = a_{1}(x,t)q_{xx}(x,t) + a_{2}(x,t)q(x,t)^{2}\overline{q(-x,t)}
\end{equation}

\noindent
respectively. We can decouple these equations into the systems

\begin{eqnarray} \label{NLS}
&& iq_{t}(x,t) + f(x,t)q_{xx}(x,t) + g(x,t)q^{2}(x,t)r(x,t) + v(x,t)q(x,t) + i\gamma(x,t)q(x,t) = 0 \\
&& -ir_{t}(x,t) + f(x,t)r_{xx}(x,t) + g(x,t)r^{2}(x,t)q(x,t) + v(x,t)r(x,t) - i\gamma(x,t)r(x,t) = 0
\end{eqnarray}

\noindent
and

\begin{eqnarray} \label{PTNLS}
iq_{t}(x,t) &=& a_{1}(x,t)q_{xx}(x,t) + a_{2}(x,t)q^{2}(x,t)r(x,t) \\
-ir_{t}(x,t) &=& a_{1}(x,t)r_{xx}(x,t) + a_{2}(x,t)r^{2}(x,t)q(x,t)
\end{eqnarray}

respectively, where in $\eqref{NLS}$ $r(x,t) = \overline{q(x,t)}$ and in $\eqref{PTNLS}$ $r(x,t) = \overline{q(-x,t)}$. It is clear that without prescribing the relation between $r(x,t)$ and $q(x,t)$, which neither Lax pair technique presented depends on, if we let $f(x,t) = -a_{1}(x,t)$, $g(x,t) = -a_{2}(x,t)$, and $v(x,t) = \gamma(x,t) = 0$ in $\eqref{NLS}$ we obtain the system given by $\eqref{PTNLS}$. That is, we can expect that the integrability conditions for the variable coefficients in the PTNLS will be a special case of the integrability conditions for variable coefficients in the NLS.

These equations, which we shall always call the physical (or field) NLPDEs to distinguish them from the many other NLPDEs we encounter, will be Lax-integrable or S-integrable if we can find a Lax pair whose compatibility condition $\eqref{ZCC1}$
contains the appropriate equation ($\eqref{NLS}$ or $\eqref{PTNLS}$). 

One expands the Lax pair $\textbf{U}$ and $\textbf{V}$ in powers of $q$ and $r$ and their derivatives with unknown functions as coefficients. This results in a large system of coupled NLPDEs for the variable coefficient functions in $\eqref{NLS}$-$\eqref{PTNLS}$. Upon solving these (and a solution is not guaranteed, and may prove to be impossible to obtain in general for some physical NLPDEs), we simultaneously obtain the Lax pair and integrability conditions on the variable coefficients for which $\eqref{NLS}$-$\eqref{PTNLS}$ are Lax-integrable. 

The results for the cubic-NLS, for which the details previously derived by Khawaja are given in Appendix A, are given by

\begin{eqnarray}
f(x,t) &=& \frac{c_{1}(t)}{g(x,t)^{2}} \\
\gamma(x,t) &=& \frac{g_{t}(x,t)}{g(x,t)} - \frac{1}{2}\frac{\dot{c_{2}}(t)}{c_{2}(t)}
\end{eqnarray}
\begin{eqnarray}
&& fg^{3}(f_{t}(g_{t}-2g\gamma)-f_{tt}g) + f_{t}^{2}g^{4} + 2f^{3}g^{3}(gv_{xx} - g_{x}v_{x}) + f^{2}g^{2}(g(4g_{t}\gamma + g_{tt}) - 2g_{t}^{2} \nonumber \\
&& - 2g^{2}(\gamma_{t} + 2\gamma^{2})) + f^{4}(36g_{x}^{4} - 48gg_{xx}g_{x}^{2} + 10g^{2}g_{xxx}g_{x} + g^{2}(6g_{xx}^{2} - gg_{xxxx})) = 0
\end{eqnarray}

The results for the PT-symmetric NLS, for which the details are also given in Appendix A, are given by

\begin{equation}
a_{1}(x,t) = \frac{h(t)}{a_{2}(x,t)^{2}}
\end{equation}

\section{The Variable Coefficient DNLS System}

Here, we will apply the technique of the last section in exactly the same fashion to generalized DNLS equation, but will omit the details for the sake of brevity. Please note that {\it the coefficients $a_i$ in this section are totally distinct or different from those given the same symbols in the
previous section. All equations in this section are thus to be read independently of those in the previous one.}

Consider the variable coefficient DNLS given by

\begin{equation}
iq_{t}(x,t) = a_{1}(x,t)q_{xx}(x,t) + ia_{2}(x,t)(q(x,t)^{2}\overline{q(x,t))_{x}}
\end{equation}

As with the previous section it will be advantageous to decouple this system into the following system

\begin{eqnarray} \label{DNLS}
iq_{t}(x,t) + a_{1}(x,t)q_{xx}(x,t) + ia_{2}(x,t)(q^{2}(x,t)r(x,t))_{x} = 0 \\
-ir_{t}(x,t) + a_{1}(x,t)r_{xx}(x,t) - ia_{2}(x,t)(r^{2}(x,t)q(x,t))_{x} = 0
\end{eqnarray} 

where $r(x,t) = \overline{q(x,t)}$. As before, we consider the variable-coefficient DNLS equation to be integrable if we can find a Lax pair which satisfies $\eqref{ZCC1}$. In the method given in (cite Khwaja) one expands the Lax pair $\textbf{U}$ and $\textbf{V}$ in powers of $q$ and $r$ and their derivatives with unknown function coefficients and require $\eqref{ZCC1}$ to be equivalent to the nonlinear system. This results in a system of coupled PDEs for the unknown coefficients for which upon solving we simultaneously obtain the Lax pair and integrability conditions on the $a_{i}$. The results, for which the details are given in Appendix B, are as follows

\begin{eqnarray}
a_{1}(x,t) &=& F_{4}(t)F_{2}(x)(c_{1} + c_{2}x) - c_{1}F_{4}(t)F_{2}(x)\int{\frac{x \ dx}{F_{2}(x)}} + c_{1}xF_{4}(t)F_{2}(x)\int{\frac{dx}{F_{2}(x)}} \\
a_{2}(x,t) &=& F_{2}(x)F_{3}(t)
\end{eqnarray}

where $F_{1-4}$ are arbitrary functions in their respective variables and $c_{1,2}$ are arbitrary constants.

Having considered these two examples of various NLS-type equations, we shall now proceed to consider
whether these results may be recovered in a more algorithmic manner. As discussed in Section 1, {\it it would be
advantageous if they could be obtained without: a. requiring to know the form of the Lax Pair for the corresponding
constant-coefficient Lax-integrable equation, and b. requiring to generalize this constant-coefficient Lax Pair
by guesswork.} Towards that end, we now proceed to consider how this may be accomplished by generalizing and extending
the Estabrook-Wahlquist technique to Lax-integrable systems with variable coefficients.

\section{The Extended Estabrook-Wahlquist Technique}

In the standard Estabrook-Wahlquist method one begins with a constant coefficient NLPDE and assumes an implicit dependence on $u(x,t)$ and its partial derivatives of the spatial and time evolution matrices ($\mathbb{F},\mathbb{G}$) involved in the linear scattering problem 
\[ \psi_{x} = \mathbb{F}\psi, \ \ \ \psi_{t} = \mathbb{G}\psi \]
The evolution matrices $\mathbb{F}$ and $\mathbb{G}$ are connected via a zero-curvature condition (independence of path in spatial and time evolution) derived by mandating $\psi_{xt} = \psi_{tx}$. That is, it requires
\[ \mathbb{F}_{t} - \mathbb{G}_{x} + [\mathbb{F},\mathbb{G}] = 0 \]
provided $u(x,t)$ satisfies the NLPDE. 

Considering the forms 
\[ \mathbb{F} = \mathbb{F}(q,r,q_{x},r_{x},q_{t},r_{t},\ldots,q_{m_{1}x,n_{1}t},r_{m_{2}x,n_{2}t}) \]

and
 
\[ \mathbb{G} = \mathbb{G}(q,r,q_{x},,r_{x},q_{t},r_{t},\ldots,q_{k_{1}x,j_{1}t},r_{k_{2}x,j_{2}t}) \]

for the space and time evolution matrices where $u_{px,qt} = \frac{\partial^{p+q}u}{\partial x^{p}\partial t^{q}}$ we see that this condition is equivalent to

\begin{eqnarray*}
&& \sum_{m_{1},n_{1}}{\mathbb{F}_{q_{m_{1}x,n_{1}t}}q_{m_{1}x,(n_{1}+1)t}} + \sum_{m_{2},n_{2}}{\mathbb{F}_{r_{m_{2}x,n_{2}t}}r_{m_{2}x,(n_{2}+1)t}} - \sum_{j_{1},k_{1}}{\mathbb{G}_{q_{j_{1}x,k_{1}t}}q_{(j_{1}+1)x,k_{1}t}} \\
&& - \sum_{j_{2},k_{2}}{\mathbb{G}_{r_{j_{2}x,k_{2}t}}r_{(j_{2}+1)x,k_{2}t}} + [\mathbb{F},\mathbb{G}] = 0
\end{eqnarray*}

\noindent
From here there is often a systematic approach\cite{EW1}-\cite{EW4} to determining the form for $\mathbb{F}$ and $\mathbb{G}$ which is outlined in \cite{EW3} and will be utilized in the examples to follow. 

Typically a valid choice for dependence on $q(x,t)$, $r(x,t)$ and their partial derivatives is to take $\mathbb{F}$ to depend on all terms in the NLPDE for which there is a partial time derivative present. Similarly we may take $\mathbb{G}$ to depend on all terms for which there is a partial space derivative present. For example, given the standard NLS,

\begin{eqnarray*}
iq_{t} &=& q_{xx} + 2q^{2}r \\
ir_{t} &=& -r_{xx} - 2qr^{2}
\end{eqnarray*}

\noindent
one would consider $\mathbb{F} = \mathbb{F}(q,r)$ and $\mathbb{G} = \mathbb{G}(q,r,q_{x},r_{x})$. Imposing compatibility allows one to determine the explicit form of $\mathbb{F}$ and $\mathbb{G}$ in a very algorithmic way. Additionally the compatibility condition induces a set of constraints on the coefficient matrices in $\mathbb{F}$ and $\mathbb{G}$. These coefficient matrices subject to the constraints generate a finite dimensional matrix Lie algebra.

In the extended Estabrook-Wahlquist method we allow for $\mathbb{F}$ and $\mathbb{G}$ to be functions of $t$, $x$, $q$, $r$ and the partial derivatives of $q$ and $r$. Although the details change, the general procedure will remain essentially the same. We will begin by equating the coefficient of the highest partial derivative of the unknown function(s) to zero and work our way down until we have eliminated all partial derivatives of the unknown function(s). 

{\it This typically results in a large partial differential equation (in the standard Estabrook-Wahlquist method, this is an algebraic equation) which can be solved by equating the coefficients of the different powers of the unknown function(s) to zero.} This final step induces a set of constraints on the coefficient matrices in $\mathbb{F}$ and $\mathbb{G}$. {\it Another big difference which we will see in the examples comes in the final and, arguably, the hardest step. In the standard Estabrook-Wahlquist method the final step involves finding explicit forms for the set of coefficient matrices such that they satisfy the contraints derived in the procedure. Note these constraints are nothing more than a system of algebraic matrix equations. In the extended Estabrook-Wahlquist method these constraints will be in the form of matrix partial differential equations which can be used to derive an integrability condition on the coefficients in the NLPDE.}

As we are now letting $\mathbb{F}$ and $\mathbb{G}$ have explicit dependence on $x$ and $t$ and for notational clarity, it will be more convenient to consider the following version of the zero-curvature condition

\begin{equation}\label{ZCC}
\mbox{D}_{t}\mathbb{F} - \mbox{D}_{x}\mathbb{G} + [\mathbb{F},\mathbb{G}] = 0
\end{equation}

\noindent
where $\mbox{D}_{t}$ and $\mbox{D}_{x}$ are the total derivative operators on time and space, respectively. Recall the definition of the total derivative

\[ \mbox{D}_{y}f(y,z,u_{1}(y,z),u_{2}(y,z),\ldots,u_{n}(y,z)) = \frac{\partial f}{\partial y} + \frac{\partial f}{\partial u_{1}}\frac{\partial u_{1}}{\partial y} + \frac{\partial f}{\partial u
_{2}}\frac{\partial u_{2}}{\partial y} + \cdots + \frac{\partial f}{\partial u_{n}}\frac{\partial u_{n}}{\partial y} \]

\noindent
Thus we can write the compatibility condition as

\begin{eqnarray*}
&& \mathbb{F}_{t} + \sum_{m_{1},n_{1}}{\mathbb{F}_{q_{m_{1}x,n_{1}t}}q_{m_{1}x,(n_{1}+1)t}} + \sum_{m_{2},n_{2}}{\mathbb{F}_{r_{m_{2}x,n_{2}t}}r_{m_{2}x,(n_{2}+1)t}} - \sum_{j_{1},k_{1}}{\mathbb{G}_{q_{j_{1}x,k_{1}t}}q_{(j_{1}+1)x,k_{1}t}} \\
&& - \sum_{j_{2},k_{2}}{\mathbb{G}_{r_{j_{2}x,k_{2}t}}r_{(j_{2}+1)x,k_{2}t}} - \mathbb{G}_{x} + [\mathbb{F},\mathbb{G}] = 0
\end{eqnarray*}

\noindent
It is important to note that the subscripted $x$ and $t$ denotes the partial derivative with respect to only the $x$ and $t$ elements, respectively. That is, although $q$, $r$ and their derivatives depend on $x$ and $t$ this will not invoke use of the chain rule as they are treated as independent variables. This will become more clear in the examples of the next section. 

Note that compatibility of the time and space evolution matrices will yield a set of constraints which contain the constant coefficient constraints as a subset. In fact, taking the variable coefficients to be the appropriate constants will yield exactly the Estabrook-Wahlquist results for the constant coefficient version of the NLPDE. That is, the constraints given by the Estabrook-Wahlquist method for a constant coefficient NLPDE are always a proper subset of the constraints given by a variable-coefficient version of the NLPDE. This can easily be shown. Letting $\mathbb{F}$ and $\mathbb{G}$ not depend explicitly on $x$ and $t$ and taking the coefficients in the NLPDE to be constant the zero-curvature condition as it is written above becomes

\begin{eqnarray*}
&& \sum_{m_{1},n_{1}}{\mathbb{F}_{q_{m_{1}x,n_{1}t}}q_{m_{1}x,(n_{1}+1)t}} + \sum_{m_{2},n_{2}}{\mathbb{F}_{r_{m_{2}x,n_{2}t}}r_{m_{2}x,(n_{2}+1)t}} - \sum_{j_{1},k_{1}}{\mathbb{G}_{q_{j_{1}x,k_{1}t}}q_{(j_{1}+1)x,k_{1}t}} \\
&& - \sum_{j_{2},k_{2}}{\mathbb{G}_{r_{j_{2}x,k_{2}t}}r_{(j_{2}+1)x,k_{2}t}} + [\mathbb{F},\mathbb{G}] = 0
\end{eqnarray*}

\noindent
which is exactly the standard Estabrook-Wahlquist method. 

The conditions derived via mandating $\eqref{ZCC}$ be satisfied upon solutions of the vc-NLPDE may be used to determine conditions on the coefficient matrices and variable-coefficients (present in the NLPDE). Successful closure of these conditions is equivalent to the system being S-integrable. A major advantage to using the Estabrook-Wahlquist method that carries forward with the extension is the fact that it requires little guesswork and yields quite general results. 

In Khawaja's method\cite{Khawaja}-\cite{Lecce} an educated guess is made for the structure of the variable-coefficient NLS Lax pair based on the associated constant coefficient Lax pair. That is, Khawaja considered the matrices
\[ \mathbb{F} = U = \begin{bmatrix}
f_{1} + f_{2}q & f_{3} + f_{4}q \\
f_{5} + f_{6}r & f_{7} + f_{8}r
\end{bmatrix} \]
and
\[ \mathbb{G} = V = \begin{bmatrix}
g_{1} + g_{2}q + g_{3}q_{x} + g_{4}rq & g_{5} + g_{6}q + g_{7}q_{x} + g_{8}rq \\
g_{9} + g_{10}r + g_{11}r_{x} + g_{12}rq & g_{13} + g_{14}r + g_{15}r_{x} + g_{16}rq
\end{bmatrix} \]
where $f_{i}$ and $g_{i}$ unknown functions of $x$ and $t$ which satisfy conditions derived by mandating the zero-curvature condition be satisfied on solutions of the variable-coefficient NLS. In fact, in a previous paper Khawaja derives the associated Lax pair via a similar means where he begins with an even weaker assumption on the structure of the Lax pair. This Lax pair is omitted from the paper as it becomes clear the zero-curvature condition mandates many of the coefficients be zero. 

{\it An ideal approach would be a method which does not require knowledge of the Lax pair to an associated constant coefficient system and involves little to no guesswork. The extended Estabrook-Wahlquist does exactly this. It will be shown that the results obtained from Khawaja's method are in fact a special case of the extended Estabrook-Wahlquist method. }

We now proceed with the variable coefficient NLS and PTNLS equations as our first examples of this extended Estabrook Wahlquist method.

\section{The vc-NLS Equation Reconsidered}

Following with the procedure outlined above we choose

\[ \mathbb{F} = \mathbb{F}(x,t,q,r), \ \ \ \mathbb{G} = \mathbb{G}(x,t,r,q,r_{x},q_{x}) \]

Compatibility requires these matrices satisfy the zero-curvature conditions given by $\eqref{ZCC}$. Plugging $\mathbb{F}$ and $\mathbb{G}$ into $\eqref{ZCC}$ we have

\begin{equation} \label{Step1}
\mathbb{F}_{t} + \mathbb{F}_{r}r_{t} + \mathbb{F}_{q}q_{t} - \mathbb{G}_{x} - \mathbb{G}_{r}r_{x} - \mathbb{G}_{q}q_{x} - \mathbb{G}_{r_{x}}r_{xx} - \mathbb{G}_{q_{x}}q_{xx} + \left[\mathbb{F},\mathbb{G}\right] = 0
\end{equation}

Now requiring this be satisfied upon solutions of $\eqref{NLS}$ we follow the standard technique of eliminating $r_{t}$ and $u_{t}$ via $\eqref{NLS}$

\begin{eqnarray} \nonumber
&& \mathbb{F}_{t} - i\mathbb{F}_{r}(fr_{xx} + gr^{2}q + (\upsilon-i\gamma)r) + i\mathbb{F}_{q}(fq_{xx} + gq^{2}r + (\upsilon+i\gamma)q) \\ \label{Step2}
&& - \mathbb{G}_{x} - \mathbb{G}_{r}r_{x} - \mathbb{G}_{q}q_{x} - \mathbb{G}_{r_{x}}r_{xx} - \mathbb{G}_{q_{x}}q_{xx} + \left[\mathbb{F},\mathbb{G}\right] = 0
\end{eqnarray}

Since $\mathbb{F}$ and $\mathbb{G}$ do not depend on $q_{xx}$ or $r_{xx}$ we collect the coefficients of $q_{xx}$ and $r_{xx}$ and equate them to zero. This requires

\begin{equation}
-if\mathbb{F}_{r} - \mathbb{G}_{r_{x}} = 0, \ \ \ \ \ if\mathbb{F}_{q} - \mathbb{G}_{q_{x}} = 0
\end{equation}

Solving this linear system yields

\begin{equation}
\mathbb{G} = if(\mathbb{F}_{q}q_{x} - \mathbb{F}_{r}r_{x}) + \mathbb{K}^{0}(x,t,q,r)
\end{equation}

Plugging this into $\eqref{Step2}$ gives us

\begin{eqnarray}
&& \mathbb{F}_{t} - i\mathbb{F}_{r}(gr^{2}q + (\upsilon-i\gamma)r) + i\mathbb{F}_{q}(gq^{2}r + (\upsilon+i\gamma)q) - if_{x}(\mathbb{F}_{q}q_{x} - \mathbb{F}_{r}r_{x}) - \mathbb{K}^{0}_{q}q_{x} - \mathbb{K}^{0}_{r}r_{x} \nonumber \\ 
&& - if(\mathbb{F}_{qx}q_{x} - \mathbb{F}_{rx}r_{x}) - if(\mathbb{F}_{qq}q_{x}^{2} - \mathbb{F}_{rr}r_{x}^{2}) - \mathbb{K}^{0}_{x} + ifq_{x}\left[\mathbb{F},\mathbb{F}_{q}\right] - ifr_{x}\left[\mathbb{F},\mathbb{F}_{r}\right] + \left[\mathbb{F},\mathbb{K}^{0}\right] = 0 \label{Step3}
\end{eqnarray}

Now since $\mathbb{F}$ and $\mathbb{K}^{0}$ do not depend on $q_{x}$ and $r_{x}$ we collect the coefficients of the $q_{x}^{2}$ and $r_{x}^{2}$ and equate them to zero. This now requires

\[ if\mathbb{F}_{rr} = 0 = -if\mathbb{F}_{qq} \]

From this, it follows via simple integration

\[ \mathbb{F} = \mathbb{X}_{1}(x,t) + \mathbb{X}_{2}(x,t)q + \mathbb{X}_{3}(x,t)r + \mathbb{X}_{4}(x,t)rq \]

\noindent
where the $\mathbb{X}_{i}$ are arbitrary matrices whose elements are functions of $x$ and $t$. Plugging this into $\eqref{Step3}$ we obtain

\begin{eqnarray}
&& \mathbb{X}_{1,t} + \mathbb{X}_{2,t}q + \mathbb{X}_{3,t}r + \mathbb{X}_{4,t}rq - i(\mathbb{X}_{3} + \mathbb{X}_{4}q)(gr^{2}q + (\upsilon-i\gamma)r) + i(\mathbb{X}_{2} + \mathbb{X}_{4}r)(gq^{2}r + (\upsilon+i\gamma)q) \nonumber \\
&& - if_{x}((\mathbb{X}_{2} + \mathbb{X}_{4}r)q_{x} - (\mathbb{X}_{3} + \mathbb{X}_{4}q)r_{x}) - \mathbb{K}^{0}_{q}q_{x} - \mathbb{K}^{0}_{r}r_{x} - if((\mathbb{X}_{2,x} + \mathbb{X}_{4,x}r)q_{x} - (\mathbb{X}_{3,x} + \mathbb{X}_{4,x}q)r_{x}) \nonumber \\
&& - \mathbb{K}^{0}_{x} + if[\mathbb{X}_{1},\mathbb{X}_{2}]q_{x} + if[\mathbb{X}_{3},\mathbb{X}_{2}]rq_{x} + if[\mathbb{X}_{4},\mathbb{X}_{2}]rqq_{x} + if[\mathbb{X}_{1},\mathbb{X}_{4}]rq_{x} + if[\mathbb{X}_{2},\mathbb{X}_{4}]rqq_{x} \nonumber \\
&& + if[\mathbb{X}_{3},\mathbb{X}_{4}]r^{2}q_{x} - if[\mathbb{X}_{1},\mathbb{X}_{3}]r_{x} - if[\mathbb{X}_{2},\mathbb{X}_{3}]qr_{x} - if[\mathbb{X}_{4},\mathbb{X}_{3}]qrr_{x} - if[\mathbb{X}_{1},\mathbb{X}_{4}]qr_{x} \nonumber \\
&& - if[\mathbb{X}_{2},\mathbb{X}_{4}]q^{2}r_{x} - if[\mathbb{X}_{3},\mathbb{X}_{4}]qrr_{x} + [\mathbb{X}_{1},\mathbb{K}^{0}] + [\mathbb{X}_{2},\mathbb{K}^{0}]q + [\mathbb{X}_{3},\mathbb{K}^{0}]r + [\mathbb{X}_{4},\mathbb{K}^{0}]qr = 0 \label{Step4a}
\end{eqnarray}

Noting the antisymmetry of the commutator we can further simplify this to

\begin{eqnarray}
&& \mathbb{X}_{1,t} + \mathbb{X}_{2,t}q + \mathbb{X}_{3,t}r + \mathbb{X}_{4,t}rq - i(\mathbb{X}_{3} + \mathbb{X}_{4}q)(gr^{2}q + (\upsilon-i\gamma)r) + i(\mathbb{X}_{2} + \mathbb{X}_{4}r)(gq^{2}r + (\upsilon+i\gamma)q) \nonumber \\
&& - if_{x}((\mathbb{X}_{2} + \mathbb{X}_{4}r)q_{x} - (\mathbb{X}_{3} + \mathbb{X}_{4}q)r_{x}) - \mathbb{K}^{0}_{q}q_{x} - \mathbb{K}^{0}_{r}r_{x} - if((\mathbb{X}_{2,x} + \mathbb{X}_{4,x}r)q_{x} - (\mathbb{X}_{3,x} + \mathbb{X}_{4,x}q)r_{x})  \nonumber \\
&& - \mathbb{K}^{0}_{x} + if[\mathbb{X}_{1},\mathbb{X}_{2}]q_{x} + if[\mathbb{X}_{3},\mathbb{X}_{2}]rq_{x} + if[\mathbb{X}_{1},\mathbb{X}_{4}]rq_{x} + if[\mathbb{X}_{3},\mathbb{X}_{4}]r^{2}q_{x} - if[\mathbb{X}_{1},\mathbb{X}_{3}]r_{x} \nonumber \\
&& - if[\mathbb{X}_{2},\mathbb{X}_{3}]qr_{x} - if[\mathbb{X}_{1},\mathbb{X}_{4}]qr_{x} - if[\mathbb{X}_{2},\mathbb{X}_{4}]q^{2}r_{x} = 0 \label{Step4b}
\end{eqnarray}

As before, since the $\mathbb{X}_{i}$ and $\mathbb{K}^{0}$ do not depend on $r_{x}$ or $q_{x}$ we equate the coefficients of the $q_{x}$ and $r_{x}$ terms to zero. Thus we require

\begin{eqnarray} 
&& -if_{x}(\mathbb{X}_{2} + \mathbb{X}_{4}r) - \mathbb{K}^{0}_{q} - if(\mathbb{X}_{2,x} + \mathbb{X}_{4,x}r) + if[\mathbb{X}_{1},\mathbb{X}_{2}] + ifr[\mathbb{X}_{1},\mathbb{X}_{4}] \nonumber \\
&& - ifr[\mathbb{X}_{2},\mathbb{X}_{3}] + ifr^{2}[\mathbb{X}_{3},\mathbb{X}_{4}] = 0 \label{FindKa} \\ 
&& if_{x}(\mathbb{X}_{3} + \mathbb{X}_{4}q) - \mathbb{K}^{0}_{r} + if(\mathbb{X}_{3,x} + \mathbb{X}_{4,x}q) - if[\mathbb{X}_{1},\mathbb{X}_{3}] - ifq[\mathbb{X}_{1},\mathbb{X}_{4}] \nonumber \\
&& - ifq[\mathbb{X}_{2},\mathbb{X}_{3}] - ifq^{2}[\mathbb{X}_{2},\mathbb{X}_{4}] = 0 \label{FindKb}
\end{eqnarray}

Upon trying to integrate this system one finds that the system is in fact inconsistent. Recall that given a system of PDEs
\[ \Psi_{q} = \xi(q,r), \ \ \Psi_{r} = \eta(q,r) \]
if we are to recover $\Psi$ we must satisfy a consistency condition. That is, we must have $\xi_{r} = \Psi_{qr} = \Psi_{rq} = \eta_{q}$. In $\eqref{FindKa}$ and $\eqref{FindKb}$ we have

\begin{eqnarray}
\xi(q,r) &=& -if_{x}(\mathbb{X}_{2} + \mathbb{X}_{4}r) - if(\mathbb{X}_{2,x} + \mathbb{X}_{4,x}r) + if[\mathbb{X}_{1},\mathbb{X}_{2}] + ifr[\mathbb{X}_{1},\mathbb{X}_{4}] \nonumber \\
&& - ifr[\mathbb{X}_{2},\mathbb{X}_{3}] + ifr^{2}[\mathbb{X}_{3},\mathbb{X}_{4}] = 0 \\ 
\eta(q,r) &=& if_{x}(\mathbb{X}_{3} + \mathbb{X}_{4}q) + if(\mathbb{X}_{3,x} + \mathbb{X}_{4,x}q) - if[\mathbb{X}_{1},\mathbb{X}_{3}] - ifq[\mathbb{X}_{1},\mathbb{X}_{4}] \nonumber \\
&& - ifq[\mathbb{X}_{2},\mathbb{X}_{3}] - ifq^{2}[\mathbb{X}_{2},\mathbb{X}_{4}] = 0
\end{eqnarray}

Thus the consistency condition ($\xi_{r} = \eta_{q}$) requires

\begin{eqnarray*}
&& -if_{x}\mathbb{X}_{4} - if\mathbb{X}_{4,x} + if[\mathbb{X}_{1},\mathbb{X}_{4}] - if[\mathbb{X}_{2},\mathbb{X}_{3}] + 2if[\mathbb{X}_{3},\mathbb{X}_{4}]r = if_{x}\mathbb{X}_{4} + if\mathbb{X}_{4,x} - if[\mathbb{X}_{1},\mathbb{X}_{4}] \\
&& - if[\mathbb{X}_{2},\mathbb{X}_{3}] - 2if[\mathbb{X}_{2},\mathbb{X}_{4}]q
\end{eqnarray*}

But this means we must have

\begin{equation}
2if_{x}\mathbb{X}_{4} + 2if\mathbb{X}_{4,x} - 2if[\mathbb{X}_{1},\mathbb{X}_{4}] - 2if[\mathbb{X}_{3},\mathbb{X}_{4}](r+q) = 0
\end{equation}

One easy choice to make the system consistent, and for the purpose of demonstrating how this method can reproduce results previously obtained in the literature, is to set $\mathbb{X}_{4} = 0$. Thus the system becomes

\begin{eqnarray} 
\mathbb{K}^{0}_{q} &=& -if_{x}\mathbb{X}_{2} - if\mathbb{X}_{2,x} + if[\mathbb{X}_{1},\mathbb{X}_{2}] - ifr[\mathbb{X}_{2},\mathbb{X}_{3}] \label{Step5a} \\
\mathbb{K}^{0}_{r} &=& if_{x}\mathbb{X}_{3} + if\mathbb{X}_{3,x} - if[\mathbb{X}_{1},\mathbb{X}_{3}] - ifq[\mathbb{X}_{2},\mathbb{X}_{3}] \label{Step5b}
\end{eqnarray}

Integrating the first equation with respect to $q$ we obtain

\[ \mathbb{K}^{0} = -if_{x}\mathbb{X}_{2}q - if\mathbb{X}_{2,x}q + if[\mathbb{X}_{1},\mathbb{X}_{2}]q - if[\mathbb{X}_{2},\mathbb{X}_{3}]rq + \mathbb{K}^{*}(x,t,r) \]

Now differentiating this and mandating that it equal our previous expression for $K^{0}_{r}$ we find that $K^{*}$ must satisfy

\[ \mathbb{K}^{*}_{r} = if_{x}\mathbb{X}_{3} + if\mathbb{X}_{3,x} - if[\mathbb{X}_{1},\mathbb{X}_{3}]  \]

Integrating this expression with respect to $r$ we easily find

\[ \mathbb{K}^{*} = if_{x}\mathbb{X}_{3}r + if\mathbb{X}_{3,x}r - if[\mathbb{X}_{1},\mathbb{X}_{3}]r + \mathbb{X}_{0}(x,t) \]

Now plugging this into our previous expression for $\mathbb{K}^{0}$ we have

\begin{equation}
\mathbb{K}^{0} = if_{x}(\mathbb{X}_{3}r - \mathbb{X}_{2}q) + if(\mathbb{X}_{3,x}r - \mathbb{X}_{2,x}q) + if[\mathbb{X}_{1},\mathbb{X}_{2}]q - if[\mathbb{X}_{1},\mathbb{X}_{3}]r - if[\mathbb{X}_{2},\mathbb{X}_{3}]qr + \mathbb{X}_{0}
\end{equation}

Now plugging this into $\eqref{Step4b}$ we have

\begin{eqnarray}
&& \mathbb{X}_{1,t} + \mathbb{X}_{2,t}q + \mathbb{X}_{3,t}r - i\mathbb{X}_{3}(gr^{2}q + (\upsilon-i\gamma)r) + i\mathbb{X}_{2}(gq^{2}r + (\upsilon+i\gamma)q) - if_{xx}(\mathbb{X}_{3}r - \mathbb{X}_{2}q) \nonumber \\
&& - 2if_{x}(\mathbb{X}_{3,x}r - \mathbb{X}_{2,x}q) - if(\mathbb{X}_{3,xx}r - \mathbb{X}_{2,xx}q) - i(f[\mathbb{X}_{1},\mathbb{X}_{2}])_{x}q + i(f[\mathbb{X}_{1},\mathbb{X}_{3}])_{x}r \nonumber \\
&& - \mathbb{X}_{0,x} + if_{x}([\mathbb{X}_{1},\mathbb{X}_{3}]r - [\mathbb{X}_{1},\mathbb{X}_{2}]q) + if([\mathbb{X}_{1},\mathbb{X}_{3,x}]r - [\mathbb{X}_{1},\mathbb{X}_{2,x}]q) + if[\mathbb{X}_{1},[\mathbb{X}_{1},\mathbb{X}_{2}]]q \nonumber \\
&& - if[\mathbb{X}_{1},[\mathbb{X}_{1},\mathbb{X}_{3}]]r - if[\mathbb{X}_{1},[\mathbb{X}_{2},\mathbb{X}_{3}]]qr + [\mathbb{X}_{1},\mathbb{X}_{0}] + if_{x}[\mathbb{X}_{2},\mathbb{X}_{3}]qr + i(f[\mathbb{X}_{2},\mathbb{X}_{3}])_{x}qr \nonumber \\
&& + if[\mathbb{X}_{2},[\mathbb{X}_{1},\mathbb{X}_{2}]]q^{2} - if[\mathbb{X}_{2},[\mathbb{X}_{1},\mathbb{X}_{3}]]qr - if[\mathbb{X}_{2},[\mathbb{X}_{2},\mathbb{X}_{3}]]q^{2}r + [\mathbb{X}_{2},\mathbb{X}_{0}]q - if_{x}[\mathbb{X}_{3},\mathbb{X}_{2}]rq \nonumber \\
&& + if([\mathbb{X}_{3},\mathbb{X}_{3,x}]r^{2} - [\mathbb{X}_{3},\mathbb{X}_{2,x}]rq) + if[\mathbb{X}_{3},[\mathbb{X}_{1},\mathbb{X}_{2}]]rq - if[\mathbb{X}_{3},[\mathbb{X}_{1},\mathbb{X}_{3}]]r^{2} - if[\mathbb{X}_{3},[\mathbb{X}_{2},\mathbb{X}_{3}]]qr^{2} \nonumber \\
&& + if([\mathbb{X}_{2},\mathbb{X}_{3,x}]qr - [\mathbb{X}_{2},\mathbb{X}_{2,x}]q^{2}) + [\mathbb{X}_{3},\mathbb{X}_{0}]r = 0 \label{Step6}
\end{eqnarray}

Since the $\mathbb{X}_{i}$ are independent of $r$ and $q$ we equate the coefficients of the different powers of $r$ and $q$ to zero and thus obtain the following constraints:

\begin{eqnarray}
O(1) &:& \mathbb{X}_{1,t} - \mathbb{X}_{0,x} + \left[\mathbb{X}_{1},\mathbb{X}_{0}\right] = 0 \\
O(q) &:& \mathbb{X}_{2,t} + i\mathbb{X}_{2}(\upsilon + i\gamma) + i(f\mathbb{X}_{2})_{xx} - i(f[\mathbb{X}_{1},\mathbb{X}_{2}])_{x} - if_{x}[\mathbb{X}_{1},\mathbb{X}_{2}] - if[\mathbb{X}_{1},\mathbb{X}_{2,x}] \nonumber \\
&& + if[\mathbb{X}_{1},[\mathbb{X}_{1},\mathbb{X}_{2}]] + [\mathbb{X}_{2},\mathbb{X}_{0}] = 0 \\
O(r) &:& \mathbb{X}_{3,t} - i\mathbb{X}_{3}(\upsilon - i\gamma) - i(f\mathbb{X}_{3})_{xx} + i(f[\mathbb{X}_{1},\mathbb{X}_{3}])_{x} + if_{x}[\mathbb{X}_{1},\mathbb{X}_{3}] + if[\mathbb{X}_{1},\mathbb{X}_{3,x}] \nonumber \\
&& - if[\mathbb{X}_{1},[\mathbb{X}_{1},\mathbb{X}_{3}]] + [\mathbb{X}_{3},\mathbb{X}_{0}] = 0 \\
O(qr) &:& 2i(f[\mathbb{X}_{2},\mathbb{X}_{3}])_{x} - if[\mathbb{X}_{1},[\mathbb{X}_{2},\mathbb{X}_{3}]] + if_{x}[\mathbb{X}_{2},\mathbb{X}_{3}] - if[\mathbb{X}_{2},[\mathbb{X}_{1},\mathbb{X}_{3}]] \nonumber \\
&& + if[\mathbb{X}_{3},[\mathbb{X}_{1},\mathbb{X}_{2}]] = 0 \\
O(q^{2}) &:& if[\mathbb{X}_{2},\mathbb{X}_{2,x}] - if[\mathbb{X}_{2},[\mathbb{X}_{1},\mathbb{X}_{2}]] = 0 \label{Satisf1} \\
O(r^{2}) &:& if[\mathbb{X}_{3},\mathbb{X}_{3,x}] - if[\mathbb{X}_{3},[\mathbb{X}_{1},\mathbb{X}_{3}]] = 0 \label{Satisf2} \\
O(q^{2}r) &:& ig\mathbb{X}_{2} - if[\mathbb{X}_{2},[\mathbb{X}_{2},\mathbb{X}_{3}]] = 0 \\
O(r^{2}q) &:& ig\mathbb{X}_{3} + if[\mathbb{X}_{3},[\mathbb{X}_{2},\mathbb{X}_{3}]] = 0
\end{eqnarray}

These equations collectively determine the conditions for integrability of the system. Note that in general, as with the standard Estabrook-Wahlquist method, the solution to the above system is not unique. Provided we can find representations for the $\mathbb{X}_{i}$ and thus reduce the system down to an integrability condition on the coefficients we will obtain our Lax pair $\mathbb{F}$ and $\mathbb{G}$. We will now show how to reproduce the results given in Khawaja's paper. Let us consider Khawaja's choices, thus

\begin{equation}
\mathbb{X}_{0} = \begin{bmatrix}
g_{1} & 0 \\
0 & g_{13}
\end{bmatrix}, \ \ \mathbb{X}_{1} = \begin{bmatrix}
f_{1} & 0 \\
0 & f_{7}
\end{bmatrix}, \ \ \mathbb{X}_{2} = \begin{bmatrix}
0 & ip_{1} \\
0 & 0
\end{bmatrix}, \ \ \mathbb{X}_{3} = \begin{bmatrix}
0 & 0 \\
-ip_{2} & 0
\end{bmatrix}
\end{equation}

Plugging this into our integrability conditions yields

\begin{eqnarray} 
O(1) &:& f_{1t} - g_{1x} = 0 \label{KC1} \\ 
O(1) &:& f_{7t} - g_{13x} = 0 \label{KC2} \\ 
O(q) &:& ip_{1t} - ip_{1}(g_{1} - g_{13} - i\upsilon + \gamma) - (fp_{1})_{xx} + 2(f_{1}-f_{7})(p_{1}f)_{x} \nonumber \\ 
&& -fp_{1}(f_{1}-f_{7})^{2} + fp_{1}(f_{1}-f_{7})_{x} = 0 \label{KC3} \\ 
O(r) &:& ip_{2t} + ip_{2}(g_{1} - g_{13} - i\upsilon - \gamma) + (fp_{2})_{xx} + 2(f_{1}-f_{7})(fp_{2})_{x} \nonumber \\ 
&& + (f_{1}-f_{7})^{2}fp_{2} + fp_{2}(f_{1}-f_{7})_{x} = 0 \label{KC4} \\ 
O(qr) &:& f_{x}p_{1}p_{2} + 2(fp_{1}p_{2})_{x} = 0 \label{KC5} \\ 
O(q^{2}r) \ \mbox{and} \ O(r^{2}q) &:& g + 2fp_{1}p_{2} = 0 \label{KC6}
\end{eqnarray}

Note that $\eqref{Satisf1}$ and $\eqref{Satisf2}$ were identically satisfied. In Khawaja's paper we see $\eqref{KC1},\eqref{KC2},\eqref{KC5}$ and $\eqref{KC6}$ given exactly. To see that the other conditions are equivalent we note that in his paper he had the additional determining equations

\begin{eqnarray} 
&& (fp_{1})_{x} - fp_{1}(f_{1}-f_{7}) - g_{6} = 0 \label{EKC1} \\ 
&& (fp_{2})_{x} + fp_{2}(f_{1}-f_{7}) - g_{10} = 0 \label{EKC2} \\ 
&& g_{6}(f_{1}-f_{7}) - ip_{1}(g_{1} - g_{13} - i\upsilon + \gamma) - g_{6x} + ip_{1t} = 0 \label{EKC3} \\ 
&& g_{10}(f_{1}-f_{7}) + ip_{2}(g_{1} - g_{13} - i\upsilon - \gamma) + g_{10x} + ip_{2t} = 0 \label{EKC4}
\end{eqnarray}

We begin by solving $\eqref{EKC1}$ and $\eqref{EKC2}$ for $g_{6}$ and $g_{10}$, respectively. Now plugging $g_{6}$ into $\eqref{EKC3}$ and $g_{10}$ into $\eqref{EKC4}$ we obtain

\begin{eqnarray}
&& 2(fp_{1})_{x}(f_{1}-f_{7}) - fp_{1}(f_{1}-f_{7})^{2} - ip_{1}(g_{1} - g_{13} - i\upsilon + \gamma) - (fp_{1})_{xx} \nonumber \\
&& + fp_{1}(f_{1}-f_{7})_{x} + ip_{1t} = 0 \\
&& (fp_{2})_{x}(f_{1}-f_{7}) + fp_{2}(f_{1}-f_{7})^{2} + ip_{2}(g_{1} - g_{13} - i\upsilon - \gamma) + (fp_{2})_{xx} \nonumber \\
&& + fp_{2}(f_{1}-f_{7})_{x} + ip_{2t} = 0
\end{eqnarray}

which is exactly $\eqref{KC3}$ and $\eqref{KC4}$. The Lax pair for this system is given by

\begin{eqnarray}
F &=& \mathbb{X}_{1} + \mathbb{X}_{2}q + \mathbb{X}_{3}r \\
G &=& if(\mathbb{X}_{2}q_{x} - \mathbb{X}_{3}r_{x}) + if_{x}(\mathbb{X}_{3}r - \mathbb{X}_{2}q) + if(\mathbb{X}_{3,x}r - \mathbb{X}_{2,x}q) + if[\mathbb{X}_{1},\mathbb{X}_{2}]q - if[\mathbb{X}_{1},\mathbb{X}_{3}]r \nonumber \\
&& - if[\mathbb{X}_{2},\mathbb{X}_{3}]qr + \mathbb{X}_{0}
\end{eqnarray}

\subsection{The Variable Coefficient PT-symmetric Nonlinear Schrodinger Equation}

As the procedure for deriving the Lax pair and conditions necessary for Lax integrability are merely a special case of the standard cubic-NLS considered in the previous section we will not go through the entire procedure but rather list the main results here. The final constraints obtained in the extended Estabrook-Wahlquist technique are the following

\begin{eqnarray}
O(1) &:& \mathbb{X}_{1,t} - \mathbb{X}_{0,x} + \left[\mathbb{X}_{1},\mathbb{X}_{0}\right] = 0 \\
O(q) &:& \mathbb{X}_{2,t} - i(a_{1}\mathbb{X}_{2})_{xx} + i(a_{1}[\mathbb{X}_{1},\mathbb{X}_{2}])_{x} + ia_{1x}[\mathbb{X}_{1},\mathbb{X}_{2}] + ia_{1}[\mathbb{X}_{1},\mathbb{X}_{2,x}] \nonumber \\
&& + if[\mathbb{X}_{1},[\mathbb{X}_{1},\mathbb{X}_{2}]] + [\mathbb{X}_{2},\mathbb{X}_{0}] = 0 \\
O(r) &:& \mathbb{X}_{3,t} + i(a_{1}\mathbb{X}_{3})_{xx} - i(a_{1}[\mathbb{X}_{1},\mathbb{X}_{3}])_{x} - ia_{1x}[\mathbb{X}_{1},\mathbb{X}_{3}] - ia_{1}[\mathbb{X}_{1},\mathbb{X}_{3,x}] \nonumber \\
&& ia_{1}[\mathbb{X}_{1},[\mathbb{X}_{1},\mathbb{X}_{3}]] + [\mathbb{X}_{3},\mathbb{X}_{0}] = 0 \\
O(qr) &:& - 2i(a_{1}[\mathbb{X}_{2},\mathbb{X}_{3}])_{x} + ia_{1}[\mathbb{X}_{1},[\mathbb{X}_{2},\mathbb{X}_{3}]] - ia_{1x}[\mathbb{X}_{2},\mathbb{X}_{3}] + ia_{1}[\mathbb{X}_{2},[\mathbb{X}_{1},\mathbb{X}_{3}]] \nonumber \\
&& - ia_{1}[\mathbb{X}_{3},[\mathbb{X}_{1},\mathbb{X}_{2}]] = 0 \\
O(q^{2}) &:& - ia_{1}[\mathbb{X}_{2},\mathbb{X}_{2,x}] + ia_{1}[\mathbb{X}_{2},[\mathbb{X}_{1},\mathbb{X}_{2}]] = 0 \label{PTSatisf1} \\
O(r^{2}) &:& - ia_{1}[\mathbb{X}_{3},\mathbb{X}_{3,x}] + ia_{1}[\mathbb{X}_{3},[\mathbb{X}_{1},\mathbb{X}_{3}]] = 0 \label{PTSatisf2} \\
O(q^{2}r) &:& - ia_{2}\mathbb{X}_{2} + ia_{1}[\mathbb{X}_{2},[\mathbb{X}_{2},\mathbb{X}_{3}]] = 0 \\
O(r^{2}q) &:& - ia_{2}\mathbb{X}_{3} - ia_{1}[\mathbb{X}_{3},[\mathbb{X}_{2},\mathbb{X}_{3}]] = 0
\end{eqnarray}
 
Utilizing the same set of generators

\begin{equation}
\mathbb{X}_{0} = \begin{bmatrix}
g_{1} & 0 \\
0 & g_{13}
\end{bmatrix}, \ \ \mathbb{X}_{1} = \begin{bmatrix}
f_{1} & 0 \\
0 & f_{7}
\end{bmatrix}, \ \ \mathbb{X}_{2} = \begin{bmatrix}
0 & ip_{1} \\
0 & 0
\end{bmatrix}, \ \ \mathbb{X}_{3} = \begin{bmatrix}
0 & 0 \\
-ip_{2} & 0
\end{bmatrix}
\end{equation}

we have the following set of conditions

\begin{eqnarray} 
O(1) &:& f_{1t} - g_{1x} = 0 \label{PTKC1} \\ 
O(1) &:& f_{7t} - g_{13x} = 0 \label{PTKC2} \\ 
O(q) &:& ip_{1t} - ip_{1}(g_{1} - g_{13}) + (a_{1}p_{1})_{xx} - 2(f_{1}-f_{7})(p_{1}a_{1})_{x} \nonumber \\ 
&& + a_{1}p_{1}(f_{1}-f_{7})^{2} - a_{1}p_{1}(f_{1}-f_{7})_{x} = 0 \label{PTKC3} \\ 
O(r) &:& ip_{2t} + ip_{2}(g_{1} - g_{13}) - (a_{1}p_{2})_{xx} - 2(f_{1}-f_{7})(a_{1}p_{2})_{x} \nonumber \\ 
&& - (f_{1}-f_{7})^{2}a_{1}p_{2} - a_{1}p_{2}(f_{1}-f_{7})_{x} = 0 \label{PTKC4} \\ 
O(qr) &:& - a_{1x}p_{1}p_{2} - 2(a_{1}p_{1}p_{2})_{x} = 0 \label{PTKC5} \\ 
O(q^{2}r) \ \mbox{and} \ O(r^{2}q) &:& a_{2} + 2a_{1}p_{1}p_{2} = 0 \label{PTKC6}
\end{eqnarray}

The Lax pair for this system can thus be found to be

\begin{eqnarray}
F &=& \mathbb{X}_{1} + \mathbb{X}_{2}q + \mathbb{X}_{3}r \\
G &=& - ia_{1}(\mathbb{X}_{2}q_{x} - \mathbb{X}_{3}r_{x}) - ia_{1x}(\mathbb{X}_{3}r - \mathbb{X}_{2}q) - ia_{1}(\mathbb{X}_{3,x}r - \mathbb{X}_{2,x}q) - ia_{1}[\mathbb{X}_{1},\mathbb{X}_{2}]q + ia_{1}[\mathbb{X}_{1},\mathbb{X}_{3}]r \nonumber \\
&& + ia_{1}[\mathbb{X}_{2},\mathbb{X}_{3}]qr + \mathbb{X}_{0}
\end{eqnarray}

Next, we illustrate our generalized Estabrook-Wahlquist technique by applying it to 
the generalized DNLS equation.

\section{The vc-DNLS Equation Reconsidered}

In this section we consider treatment of the derivative NLS given by the system

\begin{subequations} 
\begin{align}
        iq_{t} + a_{1}q_{xx} + 2iqrq_{x} + ia_{2}q^{2}r_{x}&=0,  \label{DNLSa} \\
        -ir_{t} + a_{1}r_{xx} - 2ia_{2}rq_{x}r - ia_{2}r^{2}q_{x}&=0, \label{DNLSb}
\end{align}
\end{subequations}

\noindent
where $a_{1}$ and $a_{2}$ are arbitrary functions of $x$ and $t$ with the extended Estabrook-Wahlquist method. The details of this example will be similar to that of the standard NLS and PT-symmetric NLS. Following the procedure we let

\[ \mathbb{F} = \mathbb{F}(x,t,r,q) ,\ \ \ \mathbb{G} = \mathbb{G}(x,t,r,q,r_{x},q_{x}) \]

Plugging this into $\eqref{ZCC}$ we obtain

\begin{equation}
\mathbb{F}_{t} + \mathbb{F}_{q}q_{t} + \mathbb{F}_{r}r_{t} - \mathbb{G}_{x} - \mathbb{G}_{q}q_{x} - \mathbb{G}_{q_{x}}q_{xx} - \mathbb{G}_{r}r_{x} - \mathbb{G}_{r_{x}}r_{xx} + [\mathbb{F},\mathbb{G}] = 0
\end{equation}

\noindent
Now using substituting for $q_{t}$ and $r_{t}$ using $\eqref{DNLSa}$ and $\eqref{DNLSb}$ we have

\begin{eqnarray}
&& \mathbb{F}_{t} + (ia_{1}\mathbb{F}_{q}-\mathbb{G}_{q_{x}})q_{xx} - (ia_{1}\mathbb{F}_{r}+\mathbb{G}_{r_{x}})r_{xx} - \mathbb{F}_{q}(2a_{2}rqq_{x} + a_{2}q^{2}r_{x}) \nonumber \\
&& - \mathbb{F}_{r}(2a_{2}qrr_{x} + a_{2}r^{2}q_{x}) - \mathbb{G}_{x} - \mathbb{G}_{q}q_{x} - \mathbb{G}_{r}r_{x} + [\mathbb{F},\mathbb{G}] = 0 \label{DNLSP1}
\end{eqnarray}

Since $\mathbb{F}$ and $\mathbb{G}$ do not depend on $q_{xx}$ or $r_{xx}$ we can set the coefficients of the $q_{xx}$ and $r_{xx}$ terms to zero. This requires

\begin{equation}
ia_{1}\mathbb{F}_{q} - \mathbb{G}_{q_{x}} = 0, \ \mbox{and} \ ia_{1}\mathbb{F}_{r} + \mathbb{G}_{r_{x}} = 0
\end{equation}

\noindent
Solving this in the same manner as in the NLS expample we obtain

\begin{equation}
G = ia_{1}(\mathbb{F}_{q}q_{x} - \mathbb{F}_{r}r_{x}) + \mathbb{K}^{0}(x,t,r,q)
\end{equation}

\noindent
Plugging this into $\eqref{DNLSP1}$ we obtain

\begin{eqnarray}
&& \mathbb{F}_{t} - \mathbb{F}_{q}(2a_{2}rqq_{x} + a_{2}q^{2}r_{x}) - \mathbb{F}_{r}(2a_{2}qrr_{x} + a_{2}r^{2}q_{x}) - i(a_{1}\mathbb{F}_{q})_{x}q_{x} + i(a_{1}\mathbb{F}_{r})_{x}r_{x} - \mathbb{K}^{0}_{x} \nonumber \\
&& + ia_{1}\mathbb{F}_{rr}r_{x}^{2} - \mathbb{K}^{0}_{r}r_{x} - ia_{1}\mathbb{F}_{qq}q_{x}^{2} - \mathbb{K}^{0}_{q}q_{x} + ia_{1}[\mathbb{F},\mathbb{F}_{q}]q_{x} - ia_{1}[\mathbb{F},\mathbb{F}_{r}]r_{x} + [\mathbb{F},\mathbb{K}^{0}] = 0 \label{DNLSP2}
\end{eqnarray}

Now since $\mathbb{F}$ and $\mathbb{K}^{0}$ do not depend on $q_{x}$ or $r_{x}$ we can set the coefficients of the different powers of $r_{x}$ and $q_{x}$ to zero. Thus, setting the coefficients of the $q_{x}^{2}$ and $r_{x}^{2}$ terms to zero we have

\begin{equation}
-ia_{1}\mathbb{F}_{qq} = ia_{1}\mathbb{F}_{rr} = 0
\end{equation}

\noindent
from which it follows $\mathbb{F} = \mathbb{X}_{1}(x,t) + \mathbb{X}_{2}(x,t)r + \mathbb{X}_{3}(x,t)q + \mathbb{X}_{4}(x,t)qr$. Now setting the coefficients of the $q_{x}$ and $r_{x}$ terms to zero we have

\begin{eqnarray}
&& -a_{2}\mathbb{X}_{3}q^{2} - a_{2}\mathbb{X}_{4}q^{2}r + i(a_{1}\mathbb{X}_{2})_{x} + i(a_{1}\mathbb{X}_{4})_{x}q - \mathbb{K}^{0}_{r} - ia_{1}[\mathbb{X}_{1},\mathbb{X}_{2}] - ia_{1}[\mathbb{X}_{1},\mathbb{X}_{4}]q \nonumber \\
&& - ia_{1}[\mathbb{X}_{3},\mathbb{X}_{2}]q - ia_{1}[\mathbb{X}_{3},\mathbb{X}_{4}]q^{2} - 2a_{2}\mathbb{X}_{2}rq - 2a_{2}\mathbb{X}_{4}q^{2}r = 0 \label{XIEQ2} \\
&& -a_{2}\mathbb{X}_{2}r^{2} - a_{2}\mathbb{X}_{4}r^{2}q - i(a_{1}\mathbb{X}_{3})_{x} - i(a_{1}\mathbb{X}_{4})_{x}r - \mathbb{K}^{0}_{q} + ia_{1}[\mathbb{X}_{1},\mathbb{X}_{3}] + ia_{1}[\mathbb{X}_{1},\mathbb{X}_{4}]r \nonumber \\
&& + ia_{1}[\mathbb{X}_{2},\mathbb{X}_{3}]r + ia_{1}[\mathbb{X}_{2},\mathbb{X}_{4}]r^{2} - 2a_{2}\mathbb{X}_{3}rq - 2a_{2}\mathbb{X}_{4}r^{2}q = 0 \label{ETA2}
\end{eqnarray}

In much the same way as for the NLS we denote the left-hand side of $\eqref{XIEQ2}$ as $\xi(r,q)$ and the left-hand side of $\eqref{ETA2}$ as $\eta(r,q)$. For recovery of $\mathbb{K}^{0}$ we require that $\xi_{q} = \eta_{r}$. Thus, computing $\xi_{q}$ and $\eta_{r}$ we find

\begin{eqnarray}
\xi_{q} &=& -2a_{2}\mathbb{X}_{3}q - 2a_{2}\mathbb{X}_{4}qr + i(a_{1}\mathbb{X}_{4})_{x} - ia_{1}[\mathbb{X}_{1},\mathbb{X}_{4}] - ia_{1}[\mathbb{X}_{3},\mathbb{X}_{2}] - 2ia_{1}[\mathbb{X}_{3},\mathbb{X}_{4}]q \nonumber \\
&& - 2a_{2}\mathbb{X}_{2}r - 4a_{2}\mathbb{X}_{4}qr \\
\eta_{r} &=&  -2a_{2}\mathbb{X}_{2}r - 2a_{2}\mathbb{X}_{4}rq - i(a_{1}\mathbb{X}_{4})_{x} + ia_{1}[\mathbb{X}_{1},\mathbb{X}_{4}] + ia_{1}[\mathbb{X}_{2},\mathbb{X}_{3}] + 2ia_{1}[\mathbb{X}_{2},\mathbb{X}_{4}]r \nonumber \\
&& - 2a_{2}\mathbb{X}_{3}q - 4a_{2}\mathbb{X}_{4}rq
\end{eqnarray}

\noindent
from which it follows that we must have

\begin{equation}
2i(a_{1}\mathbb{X}_{4})_{x} - 2ia_{1}[\mathbb{X}_{1},\mathbb{X}_{4}] - 2ia_{1}[\mathbb{X}_{3},\mathbb{X}_{4}]q - 2ia_{1}[\mathbb{X}_{2},\mathbb{X}_{4}]r = 0
\end{equation}

Since the $\mathbb{X}_{i}$ do not depend on $q$ or $r$ this previous condition requires

\begin{eqnarray}
&& 2i(a_{1}\mathbb{X}_{4})_{x} - 2ia_{1}[\mathbb{X}_{1},\mathbb{X}_{4}] = 0 \\
&& - 2ia_{1}[\mathbb{X}_{3},\mathbb{X}_{4}] = 0 \\
&& - 2ia_{1}[\mathbb{X}_{2},\mathbb{X}_{4}] = 0
\end{eqnarray}

As with the standard NLS we will take $\mathbb{X}_{4} = 0$ in order to simplify computations. Therefore

\begin{eqnarray}
\mathbb{K}^{0}_{q} &=& -a_{2}\mathbb{X}_{2}r^{2} - i(a_{1}\mathbb{X}_{3})_{x} + ia_{1}[\mathbb{X}_{1},\mathbb{X}_{3}] + ia_{1}[\mathbb{X}_{2},\mathbb{X}_{3}]r - 2a_{2}\mathbb{X}_{3}rq \\
\mathbb{K}^{0}_{r} &=& -a_{2}\mathbb{X}_{3}q^{2} + i(a_{1}\mathbb{X}_{2})_{x} - ia_{1}[\mathbb{X}_{1},\mathbb{X}_{2}] - ia_{1}[\mathbb{X}_{3},\mathbb{X}_{2}]q - 2a_{2}\mathbb{X}_{2}rq
\end{eqnarray}

Integrating the first equation with respect to $q$ yields

\[ \mathbb{K}^{0} = -a_{2}\mathbb{X}_{2}r^{2}q - i(a_{1}\mathbb{X}_{3})_{x}q + ia_{1}[\mathbb{X}_{1},\mathbb{X}_{3}]q + ia_{1}[\mathbb{X}_{2},\mathbb{X}_{3}]rq - a_{2}\mathbb{X}_{3}q^{2}r + \mathbb{K}^{*}(x,t,r) \]

Now differentiating this equation with respect to $r$ and mandating that it equal our previous expression for $\mathbb{K}^{0}_{r}$ we find that $\mathbb{K}^{*}$ must satisfy

\[ \mathbb{K}^{*}_{r} = i(a_{1}\mathbb{X}_{2})_{x} - ia_{1}[\mathbb{X}_{1},\mathbb{X}_{2}] \]

\noindent
from which it follows

\[ \mathbb{K}^{*} = i(a_{1}\mathbb{X}_{2})_{x}r - ia_{1}[\mathbb{X}_{1},\mathbb{X}_{2}]r + \mathbb{X}_{0}(x,t) \]

\noindent
and thus

\begin{eqnarray}
\mathbb{K}^{0} &=& i(a_{1}\mathbb{X}_{2})_{x}r - i(a_{1}\mathbb{X}_{3})_{x}q - ia_{1}[\mathbb{X}_{1},\mathbb{X}_{2}]r + ia_{1}[\mathbb{X}_{1},\mathbb{X}_{3}]q + ia_{1}[\mathbb{X}_{2},\mathbb{X}_{3}]rq - a_{2}\mathbb{X}_{2}r^{2}q \nonumber \\
&& - a_{2}\mathbb{X}_{3}q^{2}r + \mathbb{X}_{0}(x,t)
\end{eqnarray}

Now plugging this and our expression for $\mathbb{F}$ into $\eqref{DNLSP2}$ we get

\begin{eqnarray}
&& \mathbb{X}_{1,t} + \mathbb{X}_{2,t}r + \mathbb{X}_{3,t}q - i(a_{1}\mathbb{X}_{2})_{xx}r + i(a_{1}\mathbb{X}_{3})_{xx}q + i(a_{1}[\mathbb{X}_{1},\mathbb{X}_{2}])_{x}r - i(a_{1}[\mathbb{X}_{1},\mathbb{X}_{3}])_{x}q \nonumber \\
&& - i(a_{1}[\mathbb{X}_{2},\mathbb{X}_{3}])_{x}rq + (a_{2}\mathbb{X}_{2})_{x}r^{2}q + (a_{2}\mathbb{X}_{3})_{x}q^{2}r - \mathbb{X}_{0,x} + i[\mathbb{X}_{1},(a_{1}\mathbb{X}_{2})_{x}]r - [\mathbb{X}_{1},(a_{1}\mathbb{X}_{3})_{x}]q \nonumber \\
&& - ia_{1}[\mathbb{X}_{1},[\mathbb{X}_{1},\mathbb{X}_{2}]]r + ia_{1}[\mathbb{X}_{1},[\mathbb{X}_{1},\mathbb{X}_{3}]]q + ia_{1}[\mathbb{X}_{1},[\mathbb{X}_{2},\mathbb{X}_{3}]]rq - a_{2}[\mathbb{X}_{1},\mathbb{X}_{2}]r^{2}q - a_{2}[\mathbb{X}_{1},\mathbb{X}_{3}]q^{2}r \nonumber \\
&& + [\mathbb{X}_{1},\mathbb{X}_{0}] + i[\mathbb{X}_{2},(a_{1}\mathbb{X}_{2})_{x}]r^{2} - i[\mathbb{X}_{2},(a_{1}\mathbb{X}_{3})_{x}]rq - ia_{1}[\mathbb{X}_{2},[\mathbb{X}_{1},\mathbb{X}_{2}]]r^{2} + ia_{1}[\mathbb{X}_{2},[\mathbb{X}_{1},\mathbb{X}_{3}]]rq \nonumber \\
&& + ia_{1}[\mathbb{X}_{2},[\mathbb{X}_{2},\mathbb{X}_{3}]]r^{2}q + [\mathbb{X}_{2},\mathbb{X}_{0}]r + i[\mathbb{X}_{3},(a_{1}\mathbb{X}_{2})_{x}]rq - i[\mathbb{X}_{3},(a_{1}\mathbb{X}_{3})_{x}]q^{2} - ia_{1}[\mathbb{X}_{3},[\mathbb{X}_{1},\mathbb{X}_{2}]]rq \nonumber \\
&& + ia_{1}[\mathbb{X}_{3},[\mathbb{X}_{1},\mathbb{X}_{3}]]q^{2} + ia_{1}[\mathbb{X}_{3},[\mathbb{X}_{2},\mathbb{X}_{3}]]q^{2}r + [\mathbb{X}_{3},\mathbb{X}_{0}]q = 0
\end{eqnarray}

Since the $\mathbb{X}_{i}$ are independent of $r$ and $q$ we equate the coefficients of the different powers of $r$ and $q$ to zero and thus obtain the following constraints:

\begin{eqnarray}
O(1) &:& \mathbb{X}_{1,t} - \mathbb{X}_{0,x} + [\mathbb{X}_{1},\mathbb{X}_{0}] = 0 \label{DNLSCOND1} \\
O(q) &:& \mathbb{X}_{3,t} + i(a_{1}\mathbb{X}_{3})_{xx} - i(a_{1}[\mathbb{X}_{1},\mathbb{X}_{3}])_{x} - i[\mathbb{X}_{1},(a_{1}\mathbb{X}_{3})_{x}] + ia_{1}[\mathbb{X}_{1},[\mathbb{X}_{1},\mathbb{X}_{3}]] \nonumber \\
&& + [\mathbb{X}_{3},\mathbb{X}_{0}] = 0 \label{DNLSCOND2} \\
O(r) &:& \mathbb{X}_{2,t} - i(a_{1}\mathbb{X}_{2})_{xx} + i(a_{1}[\mathbb{X}_{1},\mathbb{X}_{2}])_{x} + i[\mathbb{X}_{1},(a_{1}\mathbb{X}_{2})_{x}] - ia_{1}[\mathbb{X}_{1},[\mathbb{X}_{1},\mathbb{X}_{2}]] \nonumber \\
&& + [\mathbb{X}_{2},\mathbb{X}_{0}] = 0 \label{DNLSCOND3} \\
O(rq) &:& - (a_{1}[\mathbb{X}_{2},\mathbb{X}_{3}])_{x} + a_{1}[\mathbb{X}_{1},[\mathbb{X}_{2},\mathbb{X}_{3}]] - [\mathbb{X}_{2},(a_{1}\mathbb{X}_{3})_{x}] + a_{1}[\mathbb{X}_{2},[\mathbb{X}_{1},\mathbb{X}_{3}]] + [\mathbb{X}_{3},(a_{1}\mathbb{X}_{2})_{x}] \nonumber \\
&&  - a_{1}[\mathbb{X}_{3},[\mathbb{X}_{1},\mathbb{X}_{2}]] = 0 \label{DNLSCOND4} \\
O(q^{2}) &:& - [\mathbb{X}_{3},(a_{1}\mathbb{X}_{3})_{x}] + a_{1}[\mathbb{X}_{3},[\mathbb{X}_{1},\mathbb{X}_{3}]] = 0 \label{DNLSCOND5} \\
O(r^{2}) &:& [\mathbb{X}_{2},(a_{1}\mathbb{X}_{2})_{x}] - a_{1}[\mathbb{X}_{2},[\mathbb{X}_{1},\mathbb{X}_{2}]] = 0 \label{DNLSCOND6} \\
O(r^{2}q) &:& (a_{2}\mathbb{X}_{2})_{x} - a_{2}[\mathbb{X}_{1},\mathbb{X}_{2}] + ia_{1}[\mathbb{X}_{2},[\mathbb{X}_{2},\mathbb{X}_{3}]] = 0 \label{DNLSCOND7} \\
O(q^{2}r) &:& (a_{2}\mathbb{X}_{3})_{x} - a_{2}[\mathbb{X}_{1},\mathbb{X}_{3}] + ia_{1}[\mathbb{X}_{3},[\mathbb{X}_{2},\mathbb{X}_{3}]] = 0 \label{DNLSCOND8}
\end{eqnarray}

Allowing the following forms for the generators

\begin{equation}
\mathbb{X}_{0} = \begin{bmatrix}
g_{1} & g_{2} \\
g_{3} & g_{4}
\end{bmatrix}, \ \ \mathbb{X}_{1} = \begin{bmatrix}
f_{1} & 0 \\
0 & f_{2}
\end{bmatrix}, \ \ \mathbb{X}_{2} = \begin{bmatrix}
0 & f_{3} \\
0 & 0
\end{bmatrix}, \ \ \mathbb{X}_{3} = \begin{bmatrix}
0 & 0 \\
f_{4} & 0
\end{bmatrix}
\end{equation}

Note that with this choice the $\eqref{DNLSCOND5}$ and $\eqref{DNLSCOND6}$ equations are immediately satisfed. From $\eqref{DNLSCOND2}$ and $\eqref{DNLSCOND3}$ we obtain the conditions

\begin{eqnarray}
&& g_{3}f_{4} = g_{2}f_{3} = 0 \\
&& f_{4t} + i(a_{1}f_{4})_{xx} - i(a_{1}f_{4}(f_{1} - f_{2}))_{x} + (f_{2} - f_{1})(a_{1}f_{4})_{x} + ia_{1}f_{4}(f_{1} - f_{2})^{2} \nonumber \\
&& + f_{4}(g_{4} - g_{1}) = 0 \label{DNLSFINAL2} \\
&& f_{3t} - i(a_{1}f_{3})_{xx} - i(a_{1}f_{3}(f_{1} - f_{2}))_{x} + (f_{2} - f_{1})(a_{1}f_{3})_{x} - ia_{1}f_{3}(f_{1} - f_{2})^{2} \nonumber \\
&& - f_{3}(g_{4} - g_{1}) = 0 \label{DNLSFINAL3}
\end{eqnarray}

To keep $\mathbb{X}_{2}$ and $\mathbb{X}_{3}$ nonzero we force $g_{2} = g_{3} = 0$. The condition given by $\eqref{DNLSCOND4}$ becomes the single equation

\begin{equation}
(a_{1}f_{3}f_{4})_{x} + f_{3}(a_{1}f_{4})_{x} + f_{4}(a_{1}f_{3})_{x} = 0 \label{DNLSFINAL4} \\
\end{equation}

The final two conditions now yield the system

\begin{eqnarray}
&& (a_{2}f_{3})_{x} - a_{2}f_{3}(f_{2} - f_{1}) - 2ia_{1}f_{3}^{2}f_{4} = 0 \label{DNLSFINAL5} \\
&& (a_{2}f_{4})_{x} + a_{2}f_{4}(f_{2} - f_{1}) + 2ia_{1}f_{4}^{2}f_{3} = 0 \label{DNLSFINAL6}
\end{eqnarray}

At this point solution of the system given by $\eqref{DNLSCOND1}$ and $\eqref{DNLSFINAL2} - \eqref{DNLSFINAL6}$ such that the $a_{i}$ are real-valued requires either $f_{3} = 0$ or $f_{4} = 0$. Without loss of generality we choose $f_{3} = 0$ from which we obtain the new system of equations

\begin{eqnarray}
&& f_{1t} - g_{1x} = 0 \label{DNLSE1} \\
&& f_{2t} - g_{4x} = 0 \label{DNLSE2} \\
&& f_{4t} + i(a_{1}f_{4})_{xx} - i(a_{1}f_{4}(f_{1} - f_{2}))_{x} + (f_{2} - f_{1})(a_{1}f_{4})_{x} + ia_{1}f_{4}(f_{1} - f_{2})^{2} \nonumber \\
&& + f_{4}(g_{4} - g_{1}) = 0 \label{DNLSE3} \\
&& (a_{2}f_{4})_{x} + a_{2}f_{4}(f_{2} - f_{1}) = 0 \label{DNLSE4}
\end{eqnarray}

Solving $\eqref{DNLSE1},\eqref{DNLSE3}$ and $\eqref{DNLSE4}$ for $f_{1},g_{4}$ and $f_{2}$, respectively we obtain

\begin{eqnarray}
&&f_{1} = \int{g_{1x}dt} + F_{1}(x) \\
&& f_{2} = -\frac{(a_{2}f_{4})_{x}}{a_{2}f_{4}} + \int{g_{1x}dt} + F_{1}(x) \\
&& g_{4} = \frac{-ia_{2}^{2}f_{4}a_{1xx} + ia_{1}a_{2}f_{4}a_{2xx} - 2ia_{1}a_{2x}^{2}f_{4} + 2ia_{2}a_{2x}a_{1x}f_{4} - f_{4t}a_{2}^{2}}{a_{2}^{2}f_{4}} + g_{1}
\end{eqnarray}

Plugging this into $\eqref{DNLSE2}$ yields the integrability condition

\begin{eqnarray}
&& a_{2}^{3}a_{1xxx} - ia_{2t}a_{2x}a_{2} + ia_{2xt}a_{2}^{2} - 3a_{2}^{2}a_{2xx}a_{1x} - 4a_{2x}^{3}a_{1} + 5a_{1}a_{2}a_{2x}a_{2xx} + 4a_{2x}^{2}a_{2}a_{1x} \nonumber \\
&& - a_{2}^{2}a_{1}a_{2xxx} - 2a_{2x}a_{2}^{2}a_{1xx} = 0
\end{eqnarray}

Since we require that the $a_{i}$ be real we decouple this final equation into the conditions

\begin{eqnarray}
&& a_{2t}a_{2x} - a_{2xt}a_{2} = 0 \\
&& a_{2}^{3}a_{1xxx} - 3a_{2}^{2}a_{2xx}a_{1x} - 4a_{2x}^{3}a_{1} + 5a_{1}a_{2}a_{2x}a_{2xx} + 4a_{2x}^{2}a_{2}a_{1x} \nonumber \\
&& - a_{2}^{2}a_{1}a_{2xxx} - 2a_{2x}a_{2}^{2}a_{1xx} = 0
\end{eqnarray}

With the aid of MAPLE we find that the previous system is solvable with solution given by

\begin{eqnarray}
a_{1}(x,t) &=& F_{4}(t)F_{2}(x)(c_{1} + c_{2}x) - c_{1}F_{4}(t)F_{2}(x)\int{\frac{x \ dx}{F_{2}(x)}} + c_{1}xF_{4}(t)F_{2}(x)\int{\frac{dx}{F_{2}(x)}} \\
a_{2}(x,t) &=& F_{2}(x)F_{3}(t)
\end{eqnarray}

The Lax pair for this system is thus found to be

\begin{eqnarray}
F &=& \mathbb{X}_{1} + \mathbb{X}_{3}q \\
G &=& ia_{1}\mathbb{X}_{3}q_{x} - i(a_{1}\mathbb{X}_{3})_{x}q + ia_{1}[\mathbb{X}_{1},\mathbb{X}_{3}]q - a_{2}\mathbb{X}_{3}q^{2}r + \mathbb{X}_{0}
\end{eqnarray}

This completes the extended EW analysis of the vc-DNLS equation.

\section{Conclusions and Future Work}

We have used two direct methods to obtain very significantly extended Lax- or S-integrable families
of generalized NLS, PT-symmetric NLS, and DNLS equations with coefficients which may in general vary in both space and time.
Of these, the second technique which was developed here is a new, significantly extended version of the
well-known Estabrook-Wahlquist technique for Lax-integrable systems with constant coefficients.

Future work will address the derivation of additional solutions by various methods,
as well as detailed investigations of other integrability properties of these
novel integrable inhomogeneous NLPDEs such as Backlund
Transformations, conservation laws, and, if at all feasible, bi-Hamiltonian structures and Liouville integrability.



\appendix

\newpage
\setcounter{equation}{0}
\renewcommand{\theequation}{A.\arabic{equation}}

\section{Appendix: Lax integrability for vc-NLS and vc-PTNLS equations}

We would like to note before beginning this section that the following is taken from Khawaja's paper and is included for the sake of completeness. The Lax pair is expanded in powers of $q$ and $r$ and their derivatives as follows:

\begin{eqnarray}
U &=& \begin{bmatrix}
	f_{1}+f_{2}q & f_{3}+f_{4}q  \\
	 f_{5}+f_{6}r & f_{7}+f_{8}r
\end{bmatrix} \\
V &=& \begin{bmatrix}
	g_{1} + g_{2}q + g_{3}q_{x} + g_{4}qr & g_{5} + g_{6}q + g_{7}q_{x} + g_{8}qr \\
	g_{9} + g_{10}r + g_{11}r_{x} + g_{12}qr & g_{13} + g_{14}r + g_{15}r_{x} + g_{16}qr
\end{bmatrix}
\end{eqnarray}

where $f_{1-8}$ and $g_{1-16}$ are unknown functions of $x$ and $t$. The compatibility condition gives
\begin{equation}
U_{t} - V_{x} + [U,V] = \dot{0} = \begin{bmatrix}
	0 & p_{1}(x,t)F_{1}[q,r]  \\
	p_{2}(x,t)F_{2}[q,r] & 0
\end{bmatrix}
\end{equation}
where $F_{i}[q,r]$ represents the $i^{th}$ equation in $\eqref{NLS}$ and $p_{1,2}$ are arbitrary real-valus functions. It should be clear that this off-diagonal compatibility condition requires that the coefficients of the $q$ and the $r$ on the off-diagonal of $U$ be zero. Indeed upon plugging $U$ and $V$ into the compatibility condition we immediately find that compatibility requires

\begin{eqnarray}
&& f_{2} = f_{3} = f_{5} = f_{8} = g_{2} = g_{3} = g_{5} = g_{8} = g_{9} = g_{12} = g_{14} = g_{15} = 0, f_{4} = ip_{1}, f_{6} = -ip_{2}, \nonumber \\
&& g_{7} = -fp_{1}, g_{11} = -fp_{2}, g_{4} = -g_{16} = -ifp_{1}p_{2} \nonumber
\end{eqnarray}

The remaining constraints are given by

\begin{eqnarray}
&& f_{1t} - g_{1x} = 0 \label{KhaNLS1} \\
&& f_{7t} - g_{13x} = 0 \label{KhaNLS2} \\
&& 2fp_{1}p_{2} + g = 0 \label{KhaNLS3} \\
&& f_{x}p_{1} - fp_{1}(f_{1} - f_{7}) + fp_{1x} - g_{6} = 0 \label{KhaNLS4} \\
&& f_{x}p_{2} + fp_{2}(f_{1} - f_{7}) + fp_{2x} - g_{10} = 0 \label{KhaNLS5} \\
&& g_{6}(f_{1} - f_{7}) - ip_{1}(g_{1} - g_{13} - iv + \gamma) - g_{6x} + ip_{1t} = 0 \label{KhaNLS6} \\
&& g_{10}(f_{1} - f_{7}) + ip_{2}(g_{1} - g_{13} - iv - \gamma) + g_{10x} + ip_{2t} = 0 \label{KhaNLS7} \\
&& (fp_{1}p_{2})_{x} + g_{10}p_{1} + g_{6}p_{2} = 0 \label{KhaNLS8}
\end{eqnarray}

For the reduction of these constraints to the conditions given in section 2 we refer the reader to $\cite{Khawaja}$. As the system given by $\eqref{PTNLS}$ is a special case of the system $\eqref{NLS}$ with $f = -a_{1}$, $g = -a_{2}$ and $v = \gamma = 0$ we can exploit the Lax pair constraints derived above for cubic-NLS to obtain the constraints for the PT-symmetric NLS. We therefore findm utilizing the same $U$ and $V$ that compatibility under $\eqref{PTNLS}$ requires

\begin{eqnarray}
&& f_{2} = f_{3} = f_{5} = f_{8} = g_{2} = g_{3} = g_{5} = g_{8} = g_{9} = g_{12} = g_{14} = g_{15} = 0, f_{4} = ip_{1}, f_{6} = -ip_{2}, \nonumber \\
&& g_{7} = a_{1}p_{1}, g_{11} = a_{1}p_{2}, g_{4} = -g_{16} = ia_{1}p_{1}p_{2} \nonumber
\end{eqnarray}

and the remaining constaints are given by

\begin{eqnarray}
&& f_{1t} - g_{1x} = 0 \label{KhaPTNLS1} \\
&& f_{7t} - g_{13x} = 0 \label{KhaPTNLS2} \\
&& 2a_{1}p_{1}p_{2} + a_{2} = 0 \label{KhaPTNLS3} \\
&& -a_{1x}p_{1} + a_{1}p_{1}(f_{1} - f_{7}) - a_{1}p_{1x} - g_{6} = 0 \label{KhaPTNLS4} \\
&& -a_{1x}p_{2} - a_{1}p_{2}(f_{1} - f_{7}) - a_{1}p_{2x} - g_{10} = 0 \label{KhaPTNLS5} \\
&& g_{6}(f_{1} - f_{7}) - ip_{1}(g_{1} - g_{13}) - g_{6x} + ip_{1t} = 0 \label{KhaPTNLS6} \\
&& g_{10}(f_{1} - f_{7}) + ip_{2}(g_{1} - g_{13}) + g_{10x} + ip_{2t} = 0 \label{KhaPTNLS7} \\
&& -(a_{1}p_{1}p_{2})_{x} + g_{10}p_{1} + g_{6}p_{2} = 0 \label{KhaPTNLS8}
\end{eqnarray}

For derivation of the integrability conditions on the variable coefficients present in the NLS and PT-NLS we again refer the reader to \cite{Khawaja}.

\newpage
\setcounter{equation}{0}
\renewcommand{\theequation}{B.\arabic{equation}}

\section{Appendix: Lax integrability conditions for vc-DNLS Equation}

{\it As mentioned in the text, the notation and calculations here refer to the treatment of
the generalized vcDNLS equation of Section 3 ONLY.}

The Lax pair $\textbf{U}$ and $\textbf{V}$ are expanded in powers of $q$ and $r$ and their derivatives as follows

\begin{eqnarray}
\textbf{U} &=& \begin{bmatrix}
	f_{1}+f_{2}q & f_{3}+f_{4}q  \\
	 f_{5}+f_{6}r & f_{7}+f_{8}r
\end{bmatrix} \\
\textbf{V} &=& \begin{bmatrix}
	g_{1}+g_{2}q+g_{3}q_{x}+g_{4}qr & g_{5}+g_{6}q+g_{7}qr+g_{8}q_{x}+g_{9}q^{2}r  \\
	g_{10}+g_{11}r+g_{12}qr+g_{13}r_{x}+g_{14}r^{2}q & g_{15}+g_{16}r+g_{17}r_{x}+g_{18}qr
\end{bmatrix}
\end{eqnarray}
where $f_{1-8}$ and $g_{1-20}$ are unknown functions of $x$ and $t$. Note that the compatibility condition
\begin{equation}
U_{t} - V_{x} + [U,V] = \dot{0} = \begin{bmatrix}
	0 & p_{1}(x,t)F[q,r]  \\
	0 & 0
\end{bmatrix}
\end{equation}
we enforce, where $F[q,r]$ represents the first equation of $\eqref{DNLS}$ and $p_{1}(x,t)$ is unknown, is chosen out of necessity. Upon considering a more standard compatibility condtion as that considered by Khawaja,
\begin{equation}
U_{t} - V_{x} + [U,V] = \dot{0} = \begin{bmatrix}
	0 & p_{1}(x,t)F_{1}[q,r]  \\
	p_{2}(x,t)F_{2}[q,r] & 0
\end{bmatrix}
\end{equation}
we find that the only solution mandates $p_{1}$ or $p_{2}$ be zero. We chose to let $p_{2} = 0$ but it is important to note that the conditions would not change if we had set $p_{1}$ equal to zero instead. Given this modification to the zero-curvature condition we immediately require

\begin{eqnarray}
&& f_{2}=f_{5}=f_{6}=f_{8}=g_{2}=g_{3}=g_{4}=g_{7}=g_{11}=g_{12}=g_{14}=g_{16}=g_{17}=g_{18}=0, \nonumber \\
&&  g_{8}=-p_{1}a_{1}, f_{4}=ip_{1}, g_{9}=ip_{1}a_{2} \nonumber
\end{eqnarray}

The remaining contraints are then given by

\begin{eqnarray}
&& f_{1t} - g_{1x} = 0 \label{KhaDNLS1} \\
&& f_{7t} - g_{15x} = 0 \label{KhaDNLS2} \\
&& (p_{1}a_{2})_{x} + p_{1}a_{2}(f_{7} - f_{1}) = 0 \label{KhaDNLS3} \\
&& (p_{1}a_{1})_{x} + p_{1}a_{1}(f_{7} - f_{1}) - g_{6} = 0 \label{KhaDNLS4} \\
&& - g_{5x} - g_{5}(f_{7} - f_{1}) = 0 \label{KhaDNLS5} \\
&& ip_{1t} - g_{6x} + ip_{1}(g_{15} - g_{1}) - g_{6}(f_{7} - f_{1}) = 0 \label{KhaDNLS6}
\end{eqnarray}

\subsubsection{Deriving a relation between $a_{1}$ and $a_{2}$}

Solving $\eqref{KhaDNLS1}$ and $\eqref{KhaDNLS3}-\eqref{KhaDNLS6}$ for $f_{1},f_{7}, g_{6}, g_{5}$ and $g_{15}$, respectively, we obtain

\begin{eqnarray}
f_{1} &=& \int{g_{1x}dt} + F(x) \\
f_{7} &=& -\frac{(p_{1}a_{2})_{x}}{p_{1}a_{2}} + f_{1} \\
g_{5} &=& H(t)e^{\int{(f_{1} - f_{7})dx}} \\
g_{6} &=& \frac{(p_{1}a_{1})_{x}a_{2} - (p_{1}a_{2})_{x}a_{1}}{a_{2}} \\
g_{15} &=& -\frac{i}{p_{1}}\left(g_{6x} - ip_{1t} - g_{6}(f_{1} - f_{7})\right) + g_{1}
\end{eqnarray}

Plugging these results into $\eqref{KhaDNLS2}$ we obtain the constraint

\begin{eqnarray}
&& 2p_{1}a_{2}^{2}a_{1x}p_{1xx} - 2p_{1x}^{2}a_{2}^{2}a_{1x} - 2p_{1}a_{2}a_{1x}p_{1x}a_{2x} + p_{1}^{2}a_{2}^{2}a_{1xxx} + 2p_{1x}a_{2}^{2}p_{1}a_{1xx} - p_{1}^{2}a_{2}a_{1x}a_{2xx} \nonumber \\
&& + 2p_{1x}^{2}a_{2}a_{1}a_{2x} + 2p_{1}a_{1}a_{2x}^{2}p_{1x} - 2p_{1x}a_{2}a_{1}p_{1}a_{2xx} + p_{1}^{2}a_{1}a_{2x}a_{2xx} - 2p_{1}a_{1}a_{2x}p_{1xx}a_{2} \nonumber \\
&& - p_{1}^{2}a_{2}a_{1}a_{2xxx} + ip_{1}^{2}(a_{2}a_{2xt} - a_{2t}a_{2x}) = 0
\end{eqnarray}

In order to have meaningful results we must require that the $a_{i}$ are real-valued functions. Thus we can decouple the last constraint into the following equations

\begin{eqnarray}
&& 2p_{1}a_{2}^{2}a_{1x}p_{1xx} - 2p_{1x}^{2}a_{2}^{2}a_{1x} - 2p_{1}a_{2}a_{1x}p_{1x}a_{2x} + p_{1}^{2}a_{2}^{2}a_{1xxx} + 2p_{1x}a_{2}^{2}p_{1}a_{1xx} - p_{1}^{2}a_{2}a_{1x}a_{2xx} \nonumber \\
&& + 2p_{1x}^{2}a_{2}a_{1}a_{2x} + 2p_{1}a_{1}a_{2x}^{2}p_{1x} - 2p_{1x}a_{2}a_{1}p_{1}a_{2xx} + p_{1}^{2}a_{1}a_{2x}a_{2xx} - 2p_{1}a_{1}a_{2x}p_{1xx}a_{2} \nonumber \\
&& - p_{1}^{2}a_{2}a_{1}a_{2xxx} = 0 \\
&& a_{2}a_{2xt} - a_{2t}a_{2x} = 0
\end{eqnarray}

Taking $p_{1} = a_{2}$ we obtain the final constraints

\begin{eqnarray}
&& a_{2t}a_{2x} - a_{2xt}a_{2} = 0 \\
&& a_{2}^{3}a_{1xxx} - 3a_{2}^{2}a_{2xx}a_{1x} - 4a_{2x}^{3}a_{1} + 5a_{1}a_{2}a_{2x}a_{2xx} + 4a_{2x}^{2}a_{2}a_{1x} \nonumber \\
&& - a_{2}^{2}a_{1}a_{2xxx} - 2a_{2x}a_{2}^{2}a_{1xx} = 0
\end{eqnarray}

With the aid of MAPLE we find that the previous system is solvable with solution given by

\begin{eqnarray}
a_{1}(x,t) &=& F_{4}(t)F_{2}(x)(c_{1} + c_{2}x) - c_{1}F_{4}(t)F_{2}(x)\int{\frac{x \ dx}{F_{2}(x)}} + c_{1}xF_{4}(t)F_{2}(x)\int{\frac{dx}{F_{2}(x)}} \\
a_{2}(x,t) &=& F_{2}(x)F_{3}(t)
\end{eqnarray}



\begin{thebibliography}{2}

\bibitem{K1} R. Grimshaw and S. Pudjaprasetya, Stud. Appl. Math., 112 (2004) 271.

\bibitem{K2} G Ei, R. Grimshaw, and A. Kamchatov, J. Fluid Mech., 585 (2007) 213.

\bibitem{K3} G. Das Sharma and M. Sarma, Phys. Plasmas 7 (2000) 3964.

\bibitem{K4} E. G. Fan, Phys. Lett A294 (2002), 26.

\bibitem{K5} R. Grimshaw, Proc. Roy. Soc., A368 (1979) 359.

\bibitem{K6} N. Joshi, Phys. Lett., A125 (1987), 456.

\bibitem{K7} Y. C. Zhang, J. Phys. A, 39 (2006) 14353.

\bibitem{K8} E. G. Fan, Phys. Lett A375 (2011), 493.



\bibitem{K9} X. Hu and Y. Chen, J. Nonlin. Math. Phys., 19 (2012) 1250002

\bibitem{K10} S. F. Tian and H. Q. Zhang, Stud. Appl. Math., doi 10.1111/sapm.12026, 13 pages (2013).

\bibitem{K11} S. K. Suslov, arXiv: 1012.3661v3 [math-ph].

\bibitem{K12} Z. Y Sun, Y. T. Gao, Y. Liu and X. Yu, Phys.Rev E84 (2011) 026606.

\bibitem{K13} J. He and Y. Li, Stud. Appl. Math., 126 (2010), 1.



\bibitem{Khawaja} U. Al Khawaja, A comparative analysis of Painlev\'e, Lax Pair, and Similarity Transformation methods in obtaining the integrability conditions of nonlinear Schroɤinger equations, Journal of Mathematical Physics, 51(2010), 
053506.

\bibitem{EW1} F. B. Estabrook and H. D. Wahlquist, J. Math. Phys., 16 (1975) 1; ibid: 17 (1976) 1293.

\bibitem{EW2} R. K. Dodd and A. P. Fordy, Proc. Roy. Soc. London A385 (1983) 389.

\bibitem{EW3} R. K. Dodd and A. P. Fordy, J. Phys. A17 (1984) 3249.

\bibitem{EW4} D. J. Kaup, Physica D1 (1980) 391.


\bibitem{K17} M. J. Ablowitz and Z. Musslimani, Integrable Nonlocal Nonlinear Schrodinger Equation, Physical Review Letters 
110, (2013), 064105 (4 pages).

\bibitem{K15} M. J. Ablowitz and H. Segur, Solitons and the IST (SIAM, Philadelphia, 1981).

\bibitem{K16} P. G. Drazin and R. S. Johnson, Solitons: an introduction (Cambridge U. Press, Cambridge, 1989.

\bibitem{Lecce} M. Russo and S. Roy Choudhury, J. of Phys. Conf. Series 2014: 482(1): 012038. doi 10.1088/1742-6596/482/1/012038.

\bibitem{Yu} X. Yu, Y.-T. Gao, Z.-Y. Sun, and Y. Liu, N-soliton solutions, Baɣklund transformation
and Lax pair for a generalized variable-coefficient fifth-order Korteweg Vries equation,
Phys. Scr. 81 (2010), 045402
\bibitem{Yong} X. Yong, H. Wang, and Y. Zhang, Symmetry, Intregrability and Exact Solutions of a Variable-Coefficient Korteweg-de-Vries (vcKDV) Equation, International Journal of Nonlinear Science 14 (2012)





\end{thebibliography}
\end{document}